\documentclass[12pt]{article}
\usepackage{latexsym, amssymb, amscd, amsmath}
\def\rmapdown#1{\Big\downarrow
   \rlap{$\vcenter{\hbox{$\scriptstyle#1$}}$ }}

\begin{document}
\title{\textbf{
Regular components of moduli spaces \\ of stable maps 
}}
\author{GAVRIL FARKAS
}
\date{ }
\maketitle
\newtheorem{thm}{Theorem}
\newtheorem{prop}{Proposition}[section]
\section{Introduction}

The purpose of this note is to prove the existence of `nice' components of
the Hilbert scheme of curves $C\subseteq \mathbb P^1\times \mathbb P^r$ of genus $g\geq 2$ and bidegree $(k,d)$. We can also phrase our result using the Kontsevich
moduli space of stable maps to $\mathbb P^1\times \mathbb P^r$. We work over an algebraically closed field of characteristic zero.
\newline
\indent For a smooth projective variety $Y$ and a class $\beta \in H_2(Y,\mathbb Z)$, one considers the moduli stack $\overline{\mathcal{M}}_g(Y,\beta)$ of stable maps $f:C\rightarrow Y$, with $C$ a reduced connected nodal curve of genus $g$  and $f_*([C])=\beta$ (see
[FP] for the construction of these stacks). The open substack $\mathcal{M}_g(Y,\beta)$ of $\overline{\mathcal{M}}_g(Y,\beta)$ parametrizes maps from smooth curves to $Y$. By $\overline{M}_g(Y,\beta)$ we denote the coarse moduli space corresponding to the stack $\overline{\mathcal{M}}_g(Y,\beta)$ and similarly $\overline{M}_g$ is the moduli space corresponding to the stack $\overline{\mathcal{M}}_g$ of stable curves of genus $g$. We denote by $\pi:\overline{\mathcal{M}}_g(Y,\beta)\rightarrow \overline{\mathcal{M}}_g$ the
natural projection. The {\sl expected dimension } of the stack $\overline{\mathcal{M}}_g(Y,\beta)$ is 
$$\chi(g,Y,\beta)=\mbox{dim}(Y)\ (1-g)+3g-3-\beta \cdot K_Y.$$
Since in general the geometry of $\overline{\mathcal{M}}_g(Y,\beta)$ is quite
messy (e.g. existence of many components, some nonreduced and/or not of expected dimension), it is not obvious what the definition of a nice component of $\overline{\mathcal{M}}_g(Y,\beta)$ should be. Following Sernesi [Se] we introduce the following terminology:
\newline
\noindent \textbf{Definition. } A component $V$ of $\overline{\mathcal{M}}_g(Y,\beta)$ is said to be {\sl regular} if it is generically smooth and of dimension $\chi(g,Y,\beta)$. We say that $V$ has the
{\sl expected number of moduli } if 
$$\mbox{dim }\pi(V)=\mbox{min}\bigl(3g-3,\chi(g,Y,\beta)-\mbox{dim }\mbox{Aut}(Y)\bigr).$$
\indent In this paper we only construct regular components of moduli spaces of stable maps. We study the stacks $\overline{\mathcal{M}}_g(Y,\beta)$ when $Y=\mathbb P^1\times \mathbb P^r$, $r\geq 3$ and $\beta=(k,d)\in H_2(\mathbb P^1\times \mathbb P^r,\mathbb Z).$\  We denote by $\rho(g,r,d)=g-(r+1)(g-d+r)$ the {\sl Brill-Noether number} governing the existence of $\mathfrak g^r_d$'s on curves of genus $g$. Our main result is the following:

\begin{thm}
Let $g,r,d$ and $k$ be positive integers with $r\geq 3, \rho(g,r,d)<0$ and
$$ (2-\rho(g,r,d))r+2\leq k\leq (g+2)/2.$$
Then there exists a regular component of the stack of maps $\mathcal{M}_g\bigl(\mathbb P^1\times \mathbb P^r,(k,d)\bigr)$.
\end{thm}

\indent  We introduce the {\sl Brill-Noether locus} $M_{g,d}^r=\{[C]\in M_g:C \mbox{ has a }\mathfrak g^r_d\}$,
in the case $\rho(g,r,d)<0$.  
The expected codimension of $M_{g,d}^r$ inside $M_g$ is $-\rho(g,r,d)$. We view Theorem 1 as a tool in the study of the relative position of the loci $M_{g,k}^1$ and $M_{g,d}^r$ when $r\geq 3, \rho(g,1,k)<0\Leftrightarrow k<(g+2)/2$ and $\rho(g,r,d)<0$.  The stack $\mathcal{M}_g\bigl(\mathbb P^1\times \mathbb P^r,(k,d)\bigl)$ comes naturally into play when looking at the intersection in $M_g$ of the loci $M_{g,k}^1$ and $M_{g,d}^r$. In such a setting, if $V$ is a regular component of $M_g\bigl(\mathbb P^1\times \mathbb P^r,(k,d)\bigr)$, then $M_{g,k}^1$ and $M_{g,d}^r$ intersect properly along $\pi(V)$. It is very plausible that one has a similar statement to Theorem 1 when $\rho(g,r,d)\geq 0$ and/or $\rho(g,1,k)\geq 0$, but from our perspective that seems of less interest because it would be essentially a statement about linear series on the general curve of genus $g$ with no implications on the problem of understanding the geography of the Brill-Noether loci inside $M_g$.
\newline
\indent Regarding the problem of existence of regular components of $\mathcal{M}_g(Y,\beta)$, so far the spaces $\mathcal{M}_g(\mathbb P^r,d)$
have received the bulk of attention. When $r=1,2$ the problem boils down to the study of the Hurwitz scheme and of the Severi variety of plane curves which are known to be irreducible and regular. For $r\geq 3$ we have the following result of Sernesi (cf. [Se, p. 26]):
\begin{prop}
For all $g,r,d$ such that $d\geq r+1$ and 
$$-\frac{g}{r}+\frac{r+1}{r}\leq \rho(g,r,d)<0,$$
there exists a regular component $V$ of $\mathcal{M}_g(\mathbb P^r,d)$ which has the expected number of moduli. A general point of $V$ corresponds to an embedding $C\hookrightarrow \mathbb P^r$ by a complete linear system (i.e. $h^0(C,\mathcal{O}_C(1))=r+1$),  the normal bundle $N_C$ satisfies $H^1(C,N_C)=0$ and the Petri map
$$\mu_0(C):H^0(C,\mathcal{O}_C(1))\otimes H^0(C,K_C(-1))\rightarrow H^0(C,K_C)$$
is surjective.
\end{prop}
\indent A. Lopez has obtained significant improvements on the range of $g,r,d$ such that
there exists a regular component of $\mathcal{M}_g(\mathbb P^r,d)$: if $h(r)=(4r^3+8r^2-9r+3)/(r+3)$, then for all $g,r,d$ such that $-(2-6/(r+3))g+h(r)\leq \rho(g,r,d)<0$ there exists a regular component of $\mathcal{M}_g(\mathbb P^r,d)$ with the expected number of moduli (cf. [Lo]).
\newline
\indent When  $Y$ is a smooth surface, methods from [AC] can be employed to show that if $V$ is a component of $\mathcal{M}_g(Y,\beta)$ with $\mbox{dim}(V)\geq g+1$ and which contains a point $[f:C\rightarrow Y]$ with $\mbox{deg}(f)=1$ (i.e. $f$ is generically injective), then $V$ is regular. Here it is crucial that the normal sheaf $N_f$ is of rank $1$
as then the Clifford Theorem provides an easy criterion for the vanishing of $H^1(C,N_f)$, which turns out to be a sufficient criterion for regularity (see Section 2).
\newline
\noindent \textbf{Acknowledgments.} This paper is part of my thesis written at the Universiteit van Amsterdam. The help of my adviser Gerard van der Geer is gratefully acknowledged.

\section{Deformations of maps and smoothings of space curves}

We review some facts about deformations of maps and smoothings of reducible nodal curves in $\mathbb P^r$. Our references are [Ran] and [Se].
\newline
\indent
We start by describing the deformation theory of maps between complex algebraic varieties when the source is (possibly) singular and the target is smooth.  Let $f:X\rightarrow Y$ be a morphism between complex projective varieties, with $Y$ being smooth. We denote by $\mbox{Def}(X,f,Y)$ the space of first-order deformations of the map $f$ when $X$ and $Y$ are not considered fixed. The space of first-order deformations of $X$ (resp. $Y$) is denoted by $\mbox{Def}(X)$ (resp. $\mbox{Def}(Y)$). We have the standard identification $\mbox{Def}(X)=\mbox{Ext}^1(\Omega _X, \mathcal{O}_X)$.
The deformation space $\mbox{Def}(X,f,Y)$ fits in the following exact sequence:
\begin{equation}
\mbox{Hom}_{\mathcal{O}_X}(f^*\Omega _Y,\mathcal{O}_X)\longrightarrow \mbox{Def}(X,f,Y)\longrightarrow \mbox{Def}(X)\oplus \mbox{Def}(Y)\longrightarrow \mbox{Ext}^1_f(\Omega_Y,\mathcal{O}_X).
\end{equation}
The second arrow is given by the natural forgetful maps, the space $\mbox{Hom}_{\mathcal{O}_X}(f^*\Omega _Y,\mathcal{O}_X)=H^0(X, f^*T_Y)$
parametrizes first-order deformations of $f:X\rightarrow Y$ when both $X$ and $Y$ are fixed, while for $A, B$, respectively $\mathcal{O}_X$ and $\mathcal{O}_Y$-modules, $\mbox{Ext}^i_f(B,A)$ denotes the derived functor of
$\mbox{Hom}_f(B,A)=\mbox{Hom}_{\mathcal{O}_X}(f^*B,A)=\mbox{Hom}_{\mathcal{O}_Y}(B,f_*A)$. Under reasonable assumptions (trivially satisfied when $f$ is a finite map between nodal curves) one has that $\mbox{Ext}^1_f(\Omega_Y,\mathcal{O}_X)=\mbox{Ext}^1(f^*\Omega _Y, \mathcal{O}_X)$. Using (1) it follows that when $X$ is smooth and irreducible
and $Y$ is rigid (e.g. a product of projective spaces) $\mbox{Def}(X,f,Y)=H^0(X,N_f)$, where $N_f=Coker\{T_X\rightarrow f^*T_Y\}$ is the normal sheaf of the map $f$.
 \newline
\indent For a smooth variety $Y$, a class $\beta \in H_2(Y,\mathbb Z)$ and a point $[f:C\rightarrow Y]\in \mathcal{M}_g(Y,\beta)$  we have that $T_{[f]}(\overline{\mathcal{M}}_g(Y,\beta))=H^0(C,N_f)$. If moreover $\mbox{deg}(f)=1$ and $H^1(C,N_f)=0$, then every class in $H^0(C,N_f)$ is unobstructed, $f$ is an immersion (cf. [AC, Lemma 1.4]) and $\overline{\mathcal{M}}_g(Y,\beta)$ is smooth and of the expected dimension at the 
point $[f]$, that is, $[f]$ belongs to a regular component of $\overline{\mathcal{M}}_g(Y,\beta)$.
\newline
\indent Let $C\subseteq \mathbb P^r$ be a stable curve of genus $g$ and degree $d$. If $\mathcal{I}_C$ is the ideal sheaf of $C$ we denote by $N_C:=Hom(\mathcal{I}_C/\mathcal{I}_C^2,\mathcal{O}_C)$ the normal sheaf of $C$ in $\mathbb P^r$. Assume that $H^1(C,N_C)=0$ and that $h^0(C,\mathcal{O}_C(1))=r+1$, that is, $C$ is embedded by a complete linear system. The differential of the map $\pi:\overline{\mathcal{M}}_g(\mathbb P^r,d)\rightarrow \overline{\mathcal{M}}_g$ at the point $[C\hookrightarrow \mathbb P^r]$ is given by the natural map $H^0(C,N_C)\rightarrow \mbox{Ext}^1(\Omega_C,\mathcal{O}_C)$. If $\omega_C$ denotes the dualizing sheaf of $C$, then $\mbox{rk}(d\pi)_{[C\hookrightarrow \mathbb P^r]}=3g-3-\mbox{dim }\mbox{Ker}\mu_0(C)$, where
$$\mu_0(C):H^0(C,\mathcal{O}_C(1))\otimes H^0(C,\omega_C(-1))\rightarrow H^0(C,\omega_C)$$
is the Petri map. In particular $(d\pi )_{[C\hookrightarrow \mathbb P^r]}$ has rank $3g-3+\rho(g,r,d)$ if and only if $\mu_0(C)$ is surjective.
\newline
\indent In the same setting, via the standard identification $T_{[C]}(\overline{\mathcal{M}}_g)^{\vee}=H^0(C,\omega_C\otimes \Omega_C)$, the annihilator $(\mbox{Im}(d\pi)_{[C\hookrightarrow \mathbb P^r]})^{\perp}\subseteq H^0(C,\omega_C\otimes \Omega_C)$ can be naturally identified with $\mbox{Im}(\mu_1(C))$, where 
$$\mu_1(C):\mbox{Ker}\mu_0(C)\rightarrow H^0(C,\Omega_C\otimes \omega_C)$$
is the Gaussian map obtained from taking the `derivative' of $\mu_0(C)$ (cf. [CGGH, p. 163]).
\vskip 4pt
\indent
In Section 3 we will smooth curves $X\subseteq \mathbb P^r$ which are 
unions of two smooth curves $C$ and $E$ meeting quasi-transversally (i.e. having distinct tangent lines) at a finite set $\Delta$. For such a curve one has the exact sequences (cf. [Se, p. 35])
\begin{equation} 
0\longrightarrow \mathcal{O}_E(-\Delta)\longrightarrow \mathcal{O}_X\longrightarrow \mathcal{O}_C\longrightarrow 0, 
\end{equation}

\noindent and
\begin{equation} 
0\longrightarrow \Omega _E\longrightarrow \omega _X\longrightarrow \Omega _C(\Delta)\longrightarrow 0.
\end{equation}
\indent Also in Section 3 we will use an inductive procedure to construct curves $C\subseteq \mathbb P^1\times \mathbb P^r$ with $H^1(C,N_{C/\mathbb P^1\times \mathbb P^r})=0$. The induction step uses the following result (cf. [BE, Lemma 2.3]):
\begin{prop}
Let $C\subseteq \mathbb P^r$ be a smooth curve with $H^1(C,N_C)=0$. We take $r+2$ points $p_1,\ldots,p_{r+2}\in C$ in general linear
position and a smooth rational curve $E\subseteq \mathbb P^r$ of degree $r$ which meets $C$ quasi-transversally at $p_1,\ldots,p_{r+2}$. Then $X=C\cup E$ is smoothable in $\mathbb P^r$ and $H^1(X,N_X)=0$.
\end{prop}

\section{Existence of regular components of $\mathcal{M}_g\bigl(\mathbb P^1\times \mathbb P^r, (k,d)\bigr)$}
\noindent In this section we prove the existence of regular components of $\overline{\mathcal{M}}_g\bigl(\mathbb P^1\times \mathbb P^r,(k,d)\bigr)$ in the case $k\geq r+2, d\geq r\geq 3$, and $\rho(g,r,d)<0$. We achieve this by constructing smooth curves $C\subseteq \mathbb P^1\times \mathbb P^r$ of bidegree $(k,d)$ satisfying $H^1(C,N_{C/\mathbb P^1\times \mathbb P^r})=0$. \newline
\indent 
Let us fix integers $g\geq 2, d\geq r\geq 3$ and $k\geq 2$, as well as a smooth curve $C$ of genus $g$ with maps $f_1:C\rightarrow \mathbb P^1, f_2:C\rightarrow \mathbb P^r$,  such that $\mbox{deg}(f_1)=k, \mbox{ deg}(f_2(C))=d$ and $f_2$ is generically injective. Let us denote by $f :C\rightarrow \mathbb P^1\times \mathbb P^r$ the product map. As usual we denote by $G^r_d(C)$ the scheme parametrizing $\mathfrak g^r_d$'s on $C$.
\newline
\indent There is a commutative diagram of exact sequences 
$$\begin{array}{ccccccccc}
          \; & \;  & \; & \; & \; & \; \; & 0 & \; & \; \\
           \; & \;  & \; & \; & \; & \; \; & \downarrow & \; & \; \\

          \; & \;  & \; & \; & \; & \; \; & T_C & \; & \; \\

        \; & \;  & \; & \; & \; & \; \; &\downarrow  & \; & \; \\
        0 & \longrightarrow & T_C &    \longrightarrow  & f^*(T_{\mathbb P ^1
\times \mathbb P ^r}) & \longrightarrow & N_f & \longrightarrow 0 \\
        \; & \; & \rmapdown{} & \; &\rmapdown{=} & \; &\rmapdown{} & \; & \; \\
        0 & \longrightarrow & T_C\oplus T_C &     \longrightarrow & f_1^*(T_{\mathbb P ^1})\oplus f_2^*(T_{\mathbb P ^r}) & \longrightarrow &
N_{f_1}\oplus N_{f_2} & \rightarrow 0 \\
         \; & \;  & \; & \; & \; & \; \; &\downarrow  & \; & \; \\
         \; & \;  & \; & \; & \; & \; \; & 0 & \; & \; .\\ 
\end{array}$$
By taking cohomology in the last column, we see that the condition
$H^1(C,N_{f})=0$ is equivalent to 
$H^1(C,N_{f_1})=0\mbox{ } (\mbox{automatic})\mbox{, }H^1(C,N_{f_2})=0$, and 
\begin{equation}
\mbox{Im}\{\delta_1:H^0(C,N_{f_1})\rightarrow H^1(C,T_C)\}+\mbox{Im}\{\delta_2:H^0(C,N_{f_2})\rightarrow H^1(C,T_C)\}=H^1(C,T_C),
\end{equation}
where $\delta_1$ and $\delta_2$ are coboundary maps. Condition (4) is equivalent (cf. Section 2) to
\begin{equation}
(d\pi _1)_{[f_1]}\ \bigl(T_{[f_1]}(\mathcal{M}_g(\mathbb P^1,k))\bigr)+(d\pi _2)_{[f_2]}\ \bigl(T_{[f_2]}(\mathcal{M}_g(\mathbb P^r,d))\bigr)=T_{[C]}(\mathcal{M}_g),
\end{equation}
where the projections $\pi _1:\mathcal{M}_g(\mathbb P^1,k)\rightarrow \mathcal{M}_g$ and $\pi _2:\mathcal{M}_g(\mathbb P^r,d)\rightarrow \mathcal{M}_g$ are the natural forgetful maps. Slightly abusing terminology, if $C$ is a smooth curve and
$(l_1,l_2)\in G^1_k(C)\times G^r_d(C)$ is a pair of base point free linear series on $C$, we say that $(C,l_1,l_2)$ satisfies $(5)$, if $(C,f_1,f_2)$ satifies $(5)$, where $f_1$ and $f_2$ are maps associated to $l_1$ and $l_2$.
\newline

Recall that a base point free pencil $\mathfrak g^1_k$ is said to be {\sl simple} if the induced covering $f:C\rightarrow \mathbb P^1$ has a single ramification point $x$ over each branch point and moreover $e_x(f)=2$.

We prove the existence of regular components of $\mathcal{M}_g\bigl(\mathbb P^1\times \mathbb P^r, (k,d)\bigr)$  using the following inductive procedure:
\begin{prop}
Fix positive integers $g,r,d$ and $k$ with $d\geq r\geq 3, k\geq r+2$ and $\rho(g,r,d)<0$. Let us assume that $C\subseteq \mathbb P^r$ is a smooth
nondegenerate curve of degree $d$ and genus $g$, such that $h^1(C,N_C)=0, h^0(C,\mathcal{O}_C(1))=r+1$ and the Petri map $$\mu_0(C)=\mu _0(C,\mathcal{O}_C(1)):H^0(C,\mathcal{O}_C(1))\otimes H^0(C, K_C(-1))\rightarrow H^0(C,K_C)$$ is surjective. Assume furthermore that $C$ possesses a simple base point free pencil $\mathfrak g^1_k$ say $l$, such that $|\mathcal{O}_C(1)|(-l)=\emptyset$ and $(C, l,|\mathcal{O}_C(1)|)$ satisfies (5).\newline
\indent Then there exists a smooth nondegenerate curve $Y\subseteq \mathbb P^r$
with $g(Y)=g+r+1$, $\rm{deg}$$(Y)=d+r$ and a simple base point free pencil $l'\in G^1_k(Y)$,  so that $Y$
enjoys exactly the same properties:\mbox{ } $h^1(Y,N_Y)=0,\mbox{ } h^0(Y,\mathcal{O}_Y(1))=r+1$,
the Petri map $\mu_0(Y)$ is surjective, $|\mathcal{O}_Y(1)|(-l')=\emptyset$ and
$(Y,l',|\mathcal{O}_Y(1)|)$ satifies (5).

\end{prop} 
{\sl Proof. } We first construct a reducible $k$-gonal nodal curve $X\subseteq \mathbb P^r$, with $p_a(X)=g+r+1,\mbox{ deg}(X)=d+r$, having all the required properties, then we prove that $X$ can be smoothed in $\mathbb P^r$ preserving all properties we want.
\newline
\indent Let $f_1:C\rightarrow\mathbb P^1$ be the degree $k$ map corresponding to the pencil $l$. The covering $f_1$ is simple  hence the monodromy of $f_1$ is the full symmetric group. Then since $|\mathcal{O}_C(1)|(-l)=\emptyset$, we have that for a general $\lambda \in \mathbb P^1$ the fibre $f_1^{-1}(\lambda )=p_1+\cdots +p_k$ consists of $k$ distinct points in general linear position. Let $\Delta =\{p_1,\ldots,p_{r+2}\}$ be a subset of $f_1^{-1}(\lambda )$ and let $E\subseteq \mathbb P^r$ be a rational normal curve $(\mbox{deg}(E)=r)$ passing through $p_1,\ldots,p_{r+2}$. (Through any $r+3$ points in general linear position in $\mathbb P^r$, there passes a unique rational normal curve).
We set $X:=C\cup E$, with $C$ and $E$ meeting quasi-transversally at $\Delta $.
Of course $p_a(X)=g+r+1$ and $\mbox{deg}(X)=d+r$. Note that $\rho(g,r,d)=\rho(g+r+1,r,d+r)$.\newline
\indent We first prove that $[X]\in \overline{M}_{g+r+1,k}^1$ (that is, $X$ is a limit of smooth $k$-gonal curves), by constructing an admissible covering of degree $k$ having as domain a curve $X'$, stably equivalent to $X$. Let $X':=X\cup D_{r+3}\cup\ldots\cup D_{k}$, where $D_i\simeq \mathbb P^1$ and $D_i\cap X=\{p_i\}$, for $i=r+3,\ldots,k$. Take $Y:=(\mathbb P^1)_1\cup _{\lambda}(\mathbb P^1)_2$ a union of two lines identified at $\lambda$. We construct a degree $k$ admissible covering $f':X'\rightarrow Y$ as follows: take $f'_{|{C}}=f_1:C\rightarrow (\mathbb P^1)_1$, $f'_{|{E}}=f_2:E\rightarrow (\mathbb P^1)_2$ a map of degree $r+2$ sending the points $p_1,\ldots,p_{r+2}$ to $\lambda$, and finally $f'_{|{D_i}}:D_i\simeq (\mathbb P^1)_2$ isomorphisms sending $p_i$ to $\lambda$. Clearly $f'$ is an admissible covering, so $X$ which is stably equivalent to $X'$ is a $k$-gonal curve.
\newline
\indent Let us consider now the space $\overline{\mathcal{H}}_{g+r+1,k}$ of Harris-Mumford admissible coverings of degree $k$ (cf. [HM]) and denote by $\pi_1:\overline{\mathcal{H}}_{g+r+1,k}\rightarrow \overline{\mathcal{M}}_{g+r+1}$ the natural projection which sends a covering to the stable model of its source. We assume for simplicity that $\mbox{Aut}(C)=\{Id_C\}$ which implies that  $\mbox{Aut}(f')=\{Id_{X'}\}$, so $[f']$ is a smooth point of $\overline{\mathcal{H}}_{g+r+1,k}$. In the case when $C$ has nontrivial automorphisms the argument carries through without change if we replace the space of admisible coverings with the space of twisted covers of Abramovich, Corti and Vistoli (cf. [ACV]).

We compute the differential of the map $\pi_1$ at $[f']$. We have $T_{[f']}(\overline{\mathcal{H}}_{g+r+1,k})=\mbox{Def}(X',f',Y) =
\mbox{Def}(X,f,Y)$, where $f=f'_{|_{X}}:X\rightarrow Y$. The differential 
$(d\pi_1)_{[f']}$ is the forgetful map $\mbox{Def}(X,f,Y)\rightarrow \mbox{Def}(X)$ and from the sequence (2.1) we get that $\mbox{Im}(d\pi_1)_{[f']}=
u_1^{-1}(\mbox{Im }{u_2})$, where $u_1:\mbox{Def}(X)\rightarrow \mbox{Ext}^1(f^*\Omega _Y,\mathcal{O}_X)$ and $u_2:\mbox{Def}(Y)\rightarrow \mbox{Ext}^1(f^*\Omega_Y,\mathcal{O}_X)$ are the dual maps of $u_1^{\vee}:H^0(X,\omega_X\otimes f^*\Omega_Y)\rightarrow H^0(X,\omega_X\otimes \Omega_X)$ and $u_2^{\vee}:H^0(X,\omega_X\otimes f^*\Omega_Y)\rightarrow H^0(Y,\omega_Y\otimes\Omega_Y)$. Here $u_2^{\vee}$ is induced by the trace map $\mbox{tr}:f_*\omega _X\rightarrow\omega_Y.$ Starting with the exact sequence on $X$,

$$0\longrightarrow \mbox{Tors}(\omega_X\otimes \Omega_X)
\longrightarrow \omega_X\otimes \Omega_X \longrightarrow \Omega_C^{\otimes 2}(\Delta)\oplus\Omega_E^{\otimes 2}(\Delta)\longrightarrow 0,$$

\noindent we can write the following commutative diagram of sequences
$$\begin{array}{cccccc}
     0 & \; & 0 & \; & 0 \\
    \downarrow & \; & \downarrow & \; & \downarrow \\
     H^0(\mbox{Tors}(\omega_X \otimes f^*\Omega_Y)) &    \hookrightarrow  & H^0(\omega_X\otimes f^*\Omega_Y) & \twoheadrightarrow & H^0(2K_C-R_1+\Delta)\oplus H^0(2K_E-R_2+\Delta) \\
      \rmapdown{(u_1^{\vee})_{\tiny{\mbox{tors}}}} & \; & \rmapdown{u_1^{\vee}} & \; & \rmapdown{}  \\
    H^0(\mbox{Tors}(\omega_X \otimes  \Omega_X)) &     \hookrightarrow & H^0(\omega_X\otimes \Omega_X) & \twoheadrightarrow &
H^0(2K_C+\Delta)\oplus H^0(2K_E+\Delta) \\
\end{array}$$
where $R_1$ (resp. $R_2$) is the ramification divisor of the map $f_1$ (resp. $f_2$). Taking into account that $H^0(E,2K_E-R_2+\Delta)=0$ and that $H^0(Y,\omega_Y\otimes \Omega_Y)=H^0(\mbox{Tors}(\omega_Y\otimes \Omega_Y))$, we obtain that
\begin{equation} 
\mbox{Im}(d\pi _1) _{[f']}=\bigl(H^0(C,2K_C-R_1+\Delta) \oplus  \mbox{Ker}(u_2^{\vee})_{\tiny{\mbox{tors}}}\bigr)^{\perp},
\end{equation}
where $(u_2^{\vee})_{\tiny{\mbox{tors}}}:H^0(\mbox{Tors}(\omega _X\otimes f^*\Omega_Y))\rightarrow H^0(\mbox{Tors}(\omega_Y\otimes \Omega_Y))$ is the restriction of $u_2^{\vee}$. The space $\mbox{Ker}(u_2^{\vee})_{\tiny{\mbox{tors}}}$ is just a hyperplane in $H^0(\mbox{Tors}(\omega_X\otimes f^*\Omega_Y))\simeq \mathbb C^{r+2}$.
\vskip 10pt
\noindent \textbf{Intermezzo.} If we also assume that $\rho(g,1,k)<0$ and that $[C]$ is a smooth point of $M_{g,k}^1$ (which happens precisely when $\mbox{Aut}(C)=\{Id_C\}$, $C$ has exactly one $\mathfrak g^1_k$ and $\mbox{dim}|2\mathfrak g^1_k|=2$),
then we can prove that the locus $\overline{M}_{g+r+1,k}^1$ is smooth at $[X]$ as well. Indeed, since $\Delta \in C_{r+2}$ was chosen generically in a fibre of the $\mathfrak g^1_k$ on $C$, from Riemann-Roch we have that $h^0(C,2K_C-R_1+\Delta)=g-2k+3+r=\mbox{codim}(\overline{M}^1_{g+r+1,k},\overline{M}_{g+r+1})$. The fibre over $[X]$ of the map $\pi_1:\overline{\mathcal{H}}_{g+r+1,k}\rightarrow \overline{M}_{g+r+1}$ is identified with the space of degree $r+1$ maps $f_2:E\rightarrow \mathbb P^1$ such that $f_2(p_1)=\ldots =f_2(p_{r+2})=\lambda$, hence it is $r+1$ dimensional. We compute the tangent cone

$$TC_{[X]}(\overline{M}^1_{g+r+1,k})=\bigcup \{\mbox{Im}(d\pi_1)_z:z\in \pi_1^{-1}([X])\}=H^0(C,2K_C-R_1+\Delta)^{\perp},$$
which shows that $[X]$ is a smooth point of the locus $\overline{M}_{g+r+1,k}^1$.
\vskip 8pt
\indent We compute now the differential $$(d\pi_2)_{[X]}:T_{[X]}(\mbox{Hilb}_{d+r,g+r+1,r})\rightarrow  T_{[X]}(\overline{\mathcal{M}}_{g+r+1}),$$
\noindent which is the same thing as the differential at the point $[X\hookrightarrow \mathbb P^r]$ of the projection $\pi_2:\overline{\mathcal{M}}_{g+r+1}(\mathbb P^r,d+r)\rightarrow \overline{\mathcal{M}}_{g+r+1}$. We start by noticing that $X$ is smoothable in $\mathbb P^r$ and that $H^1(X,N_X)=0$ (apply Proposition 2.1). We also have that $X$ is embedded in $\mathbb P^r$ by a complete linear system, that is,
$h^0(X,\mathcal{O}_X(1))=r+1$. Indeed, on one hand, since $X$ is nondegenerate, $h^0(X,\mathcal{O}_X(1))\geq h^0(\mathbb P^r, \mathcal{O}_{\mathbb P^r(1)})=r+1$; on the other hand from the sequence (2) we have that $h^0(X,\mathcal{O}_X(1))\leq h^0(C,\mathcal{O}_C(1))=r+1$. \newline
\indent If $X$ is embedded in $\mathbb P^r$ by a complete linear system, we know (cf. Section 2) that $$\mbox{Im}(d\pi_2)_{[X]}=(\mbox{Im}\mu_1(X))^{\perp},$$
\noindent where $\mu_1(X):\mbox{Ker} \mu_0(X)\rightarrow H^0(X,\omega_X\otimes \Omega_X)$ is the `derivative' of the Petri map $\mu_0(X):H^0(X,\mathcal{O}_X(1))\otimes H^0(X,\omega_X(-1))\rightarrow H^0(X,\omega_X)$. We compute the
kernel of $\mu_0(X)$ and show that $\mu_0(X)$ is surjective in a way that resembles the proof of Proposition 2.3 in [Se].
\newline
\indent From the sequence (3) we obtain 
$H^0(X,\omega_X)=H^0(C,K_C+\Delta)$, while from (2) we have that
$H^0(X,\mathcal{O}_X(1))=H^0(E,\mathcal{O}_E(1))$ (keeping in mind that $H^0(C,\mathcal{O}_C(1)(-\Delta))=0$, as $p_1,\ldots,p_{r+2}$ are in general linear position). Finally, using (3) again, we have that $H^0(X,\omega_X(-1))=H^0(C,K_C(-1)+\Delta )$. Therefore we can write the following commutative diagram:
$$\begin{array}{ccc}
   
     H^0(C,\mathcal{O}_C(1))\otimes H^0(C,K_C(-1)) &    \stackrel{\mu_0(C)}\longrightarrow  & H^0(C,K_C) \\
      \rmapdown{} & \; & \rmapdown{}  \\
    H^0(C,\mathcal{O}_C(1))\otimes H^0(C,K_C(-1)+\Delta) &     \longrightarrow & H^0(C,K_C+\Delta) \\
    \rmapdown{=} & \; & \rmapdown{=} \\
    H^0(X,\mathcal{O}_X(1))\otimes H^0(X,\omega_X(-1)) &\stackrel{\mu_0(X)}\longrightarrow & H^0(X,\omega_X) \mbox{ }.\\ 
\end{array}$$
It follows that $\mbox{Ker}\mu_0(C)\subseteq \mbox{Ker}\mu_0(X)$. By Corollary 1.6 from [CR], our assumptions ($\mu_0(C)$ surjective and $\mbox{card}(\Delta )\geq 4$) enable us to conclude that $\mu_0(X)$ is surjective too. Then 
$\mbox{Ker}\mu_0(C)=\mbox{Ker}\mu_0(X)$ for dimension reasons, hence also $\mbox{Im}\mu_1(X)=
\mbox{Im}\mu_1(C)\subseteq H^0(C,2K_C)\subseteq H^0(X,\omega_X\otimes \Omega_X)$. We thus get that
$\mbox{Im}(d\pi_2)_{[X]}=(\mbox{Im}\mu_1(X))^{\perp}=(\mbox{Im}\mu_1(C))^{\perp}.$\newline
\indent
The assumption that $(C,f_1,f_2)$ satisfies (5) can be rewritten by passing to duals as
$$H^0(C,2K_C-R_1)^{\perp}+(\mbox{Im}\mu_1(C))^{\perp}=H^1(C,T_C)
\Longleftrightarrow H^0(C,2K_C-R_1)\cap \mbox{Im}\mu_1(C)=0.$$
\noindent Then it follows that $\mbox{Im}\mu_1(X)\cap \bigl(H^0(C,2K_C-R_1+\Delta)\oplus \mbox{Ker}((u_2^{\vee})_{\tiny{\mbox{tors}}})\bigr )=0$, which is the same thing as

\begin{equation}
(d\pi_1)_{[f']}\ \bigl(T_{[f']}(\overline{\mathcal{H}}_{g+r+1,k})\bigr)+(d\pi_2)_{[X\hookrightarrow \mathbb P^r]}
\ \bigl(T_{[X\hookrightarrow \mathbb P^r]}(\overline{\mathcal{M}}_{g+r+1}(\mathbb P^r,d+r))\bigr)=\mbox{Ext}^1(\Omega_X,\mathcal{O}_X).
\end{equation}
This means that the images of $\overline{\mathcal{H}}_{g+r+1,k}$ under
the map $\pi_1$ and of $\overline{\mathcal{M}}_{g+r+1}(\mathbb P^r,d+r)$ under
the map $\pi_2$, meet transversally at the point $[X]\in \overline{\mathcal{M}}_{g+r+1}$.
\vskip 4pt
\noindent \textbf{Claim. } The curve $X$ can be smoothed in such 
a way that the $\mathfrak g^1_k$ and the very ample $\mathfrak g^r_{d+r}$ are preserved  (while (7) is an open condition on $\overline{\mathcal{H}}_{g+r+1,k}\times \overline{\mathcal{M}}_{g+r+1}(\mathbb P^r, d+r)$).
\vskip 6pt
\noindent Indeed, the tangent directions that fail to smooth at least one node of $X$ are those in $\bigcup_{i=1}^{r+2}H^0(\mbox{Tors}_{p_i}(\omega_X\otimes \Omega_X))^{\perp}$, whereas the tangent directions that preserve both the 
$\mathfrak g^1_k$ and the $\mathfrak g^r_{d+r}$ are those in
$$\bigl((\mbox{Im}\mu_1(C)+H^0(C,2K_C-R_1+\Delta))\oplus \mbox{Ker}(u_2^{\vee})_{\tiny{\mbox{tors}}}\bigr)^{\perp}.$$ Since 
$H^0(\mbox{Tors}_{p_i}(\omega_X\otimes \Omega_X))\nsubseteq \mbox{Ker}
(u_2^{\vee})_{\tiny{\mbox{tors}}}$ for $i=1,\ldots,r+2$, by moving in a suitable direction in the tangent space at $[f']$ of 
$\pi_1^{-1}\pi_2(\overline{\mathcal{M}}_{g+r+1}(\mathbb P^r,d+r))$, we finally obtain a smooth curve $Y\subseteq \mathbb P^r$ with $g(Y)=g+r+1, \mbox{deg}(Y)=d+r$ and 
satisfying all the required properties. \hfill $\Box$

\vskip 9pt
\indent Using the previous result together with Proposition 1.1 we construct now regular components of $\mathcal{M}_g\bigl(\mathbb P^1\times \mathbb P^r, (k,d)\bigr)$.
\vskip 4pt
\noindent {\bf Theorem 1}  {\sl Let $g,r,d$ and $k$ be positive integers such that $r\geq 3$, $\rho(g,r,d)<0$ and 
$$(2-\rho(g,r,d))r+2\leq k\leq (g+2)/2.$$
Then there exists a regular component of the stack of maps $\mathcal{M}_g\bigl(\mathbb P^1\times \mathbb P^r,(k,d)\bigr)$.}
\vskip 5pt 
\noindent {\sl Proof.} All integer solutions $(g_0,d_0)$ of the equation $\rho(g_0,r,d_0)=\rho(g,r,d)$ with $g_0\leq g$ and $d_0\leq d$, are of the form
$g_0=g-a(r+1)$ and $d_0=d-ar$ with $a\geq 0$. Using our numerical assumptions, by a routine check we find that there exists $a>0$
such that $g_0=g-a(r+1)>0$, $d_0=d-ar\geq r+1$, $k\geq g_0+1$ and
$$-\frac{g_0}{r}+\frac{r+1}{r}\leq \rho(g_0,r,d_0)<0.$$
\indent By Proposition 1.1 there exists a smooth curve $C_0\subseteq \mathbb P^r$ of genus $g_0$ and degree $d_0$, with $H^1(C_0,N_{C_0/\mathbb P^r})=0$, 
$h^0(C_0,\mathcal{O}_{C_0}(1))=r+1$ and $\mu_0(C_0)$ surjective. Moreover, since $k\geq g_0+1$, there exists an open dense subset $U\subseteq \mbox{Pic}^k(C_0)$ such that for each $L_1\in U$ there exists a pencil $l_1=(L_1,V_1)\in G^1_k(C_0)$ with $V_1\in \mbox{Gr}(2,H^0(C_0,L))$, such that $l_1$ is simple and base point free (cf. [Fu, Proposition 8.1]).
\newline

\indent

We denote by $\pi_1:\mathcal{M}_{g_0}(\mathbb P^1,k_0)\rightarrow \mathcal{M}_{g_0}$ the natural projection and by $f_1:C\rightarrow \mathbb P^1$ the map corresponding to $l_1$. By Riemann-Roch we have $H^1(C_0,L_1^{\otimes2})=0$, hence using Section 2 \ $(d\pi_1)_{[f_1]}:T_{[f_1]}(\mathcal{M}_{g_0}(\mathbb P^1,k))\rightarrow T_{[C_0]}(\mathcal{M}_{g_0})$ is surjective since $\mbox{Coker}(d\pi_1)_{[f_1]}=H^1(C_0,f_1^*T_{\mathbb P^1})=0$.
It follows that $(C_0, |\mathcal{O}_{C_0}(1)|, l_1)$ satisfies (5).
\newline
\indent We claim that if $L\in U$ is general then $|\mathcal{O}_{C_0}(1)\otimes L^{\vee}|=\emptyset$. Suppose not, that is $\mathcal{O}_{C_0}(1)\otimes L^{\vee}\in W_{d_0-k}(C_0)$ for a general $L\in \mbox{Pic}^k(C_0)$. This is possible only for 
$d_0-k\geq g_0$, hence
$$r+2\leq k\leq d_0-g_0<r, \mbox{ }(\mbox{because }\rho(g_0,r,d_0)=\rho(g,r,d)<0),$$
a contradiction. Thus $(C_0, |\mathcal{O}_{C_0}(1)|, l_1)$ satisfies all conditions required by Proposition 3.1 which we can now apply $a$ times to get a smooth curve $C\subseteq \mathbb P^1\times \mathbb P^r$ of genus $g$ and bidegree $(k,d)$ such that $H^1(C,N_{C/\mathbb P^1\times \mathbb P^r})=0$. The conclusion of Theorem 1 now follows. \hfill $\Box$
\vskip 6pt
\indent
In the special case $\rho(g,r,d)=-1$ we can extend the range of possible $g,r,d$ and $k$ for which there is a regular component:

\begin{thm}
Let $g,r,d,k$ be positive integers such that $r\geq 3$, $\rho(g,r,d)=-1$ and 
$$\frac{2r^2+r+1}{r-1}\leq k\leq \frac{g+2}{2}\mbox{ }.$$
Then there exists a regular component of the stack of maps $\mathcal{M}_g\bigl(\mathbb P^1\times \mathbb P^r, (k,d)\bigr)$.
\end{thm}

\vskip 6pt
\noindent {\sl Proof. } We find a solution $\bigl(g_0=g-a(r+1),d_0=d-ar\bigr)$
of the equation $\rho(g_0,r,d_0)=\rho(g,r,d)$ with $a\in \mathbb Z_{\geq 0}$ such that $d_0\geq k+r$ and $\rho(g_0,1,k)\geq r-1$. Our numerical assumptions ensure that such an $a\geq 0$
exists. Note that in this case $k\leq g_0+1$, so we are not in the situation covered by Theorem 1.
\newline
\indent It also follows that $-\frac{g_0}{r}+\frac{r+1}{r}\leq -1=\rho(g_0,r,d_0)$ and $d_0\geq r+1$, hence by Proposition 1.1 there exists an irreducible smooth open subset $U$ of $\mathcal{M}_{g_0}(\mathbb P^r,d_0)$ of the expected dimension, such that all points of $U$ correspond to embeddings of smooth curves
$C\hookrightarrow \mathbb P^r$, with $h^1(C,N_C)=0, \mbox{ }h^0(C,
\mathcal{O}_C(1))=r+1$ and $\mu_0(C)$ surjective. 
\newline
\indent Since we are in the case $\rho(g_0,r,d_0)=-1$, a combination of results by Eisenbud, Harris and Steffen gives that the Brill-Noether locus $M_{g_0,d_0}^{r}$ is an irreducible divisor in $M_{g_0}$ (see [St, Theorem 0.2]). It follows that the natural projection
$\pi_2:U\rightarrow \mathcal{M}_{g_0,d_0}^r$ is dominant.
\newline
\indent To apply Proposition 3.1 we now find a curve $[C_0]\in M_{g_0,d_0}^{r_0}$ having a complete base point free $\mathfrak g^1_k$ such that $2\mathfrak g^1_k$ is
non-special. Then by semicontinuity we get that the general $[C]\in U$ also possesses a pencil $\mathfrak g^1_k$ with these properties. To find one particular such curve we proceed as follows: take $C_0$ a general $(r+1)$-gonal curve of genus $g_0$. These curves will have rather few 
moduli $(r+1<[(g+3)/2]$) but we still have that $[C_0]\in M_{g_0,d_0}^{r}$. Indeed, according to [CM] we can construct a $\mathfrak g^r_{d_0}=
|\mathfrak g^1_{r+1}+F|$ on $C_0$, where $F$ is an effective divisor on $C_0$ with
$h^0(C_0,F)=1$. Since $k\leq g_0$, using Corollary 2.2.3 from [CKM] we find that $C_0$ also
carries a complete base point free $\mathfrak g^1_k$, not composed with the $\mathfrak g^1_{r+1}$ computing $\mbox{gon}(C_0)$, and such that $2\mathfrak g^1_k$ is non-special. Since these are open conditions, they will hold generically along a component of $G^1_k(C_0)$. Applying semicontinuity, for a general element $[C]\in M^r_{g_0,d_0}$ (hence also for a general element 
$[C]\in U$), the variety $G^1_k(C)$ will contain a component $A$ with general
point $l\in A$ being  complete, base point free and with $2l$ non-special. 
\newline
\indent We claim that for a general $l\in A$ we have that $|\mathcal{O}_C(1)|(-l)=\emptyset$. Suppose not. Then if we denote by $V_{d_0-k}^{r-1}(|\mathcal{O}_C(1)|)$ the variety of effective divisors of degree $d_0-k$ on $C$ imposing $\leq r-1$ conditions on $|\mathcal{O}_C(1)|$, we obtain
$$\mbox{dim }V_{d_0-k}^{r-1}(|\mathcal{O}_C(1)|)\geq \mbox{dim }A\geq \rho(g_0,1,k)\geq r-1.$$
\indent Therefore $C\subseteq \mathbb P^r$ has at least $\infty^{r-1}$ $(d_0-k)$-secant $(r-2)$-planes, hence also at least $\infty^{r-1}$
$r$-secant $(r-2)$-planes (because $d_0-k\geq r).$
This last statement contradicts the Uniform Position Theorem (see [ACGH,
p. 112]), hence the general point $[C]\in U$ enjoys all properties required to make Proposition 3.1 work. \hfill $\Box$

\vskip 10pt  
\noindent \textbf{Remark. } From the proof of Theorem 2 the following question appears naturally: let us fix $g,k$ such that $g/2+1\leq k\leq g$. One knows
(cf. [ACGH]) that if $l\in G^1_k(C)$ is a complete, base point free pencil then
$\mbox{dim}\ T_l(G^1_k(C))=\rho(g,1,k)+h^1(C,2l)$. Therefore if $A$ is a component of $G^1_k(C)$ such that $\mbox{dim\ }A=\rho(g,1,k)$ and the general $l\in A$ is base point free such that $2l$ is special, then $A$ is nonreduced.  What is then the dimension of the locus
$$V_{g,k}:=\{[C]\in M_g:\mbox{ every component of }G^1_k(C)\mbox{ is nonreduced }\}?$$
\noindent A result of Coppens (cf. [Co]) says that for a curve $C$, if the
scheme $W^1_{k}(C)$ is reduced and of dimension $\rho(g,1,k)$, then the scheme $W^1_{k+1}(C)$ is reduced too and of dimension $\rho(g,1,k+1)$. It would
make then sense to determine $\mbox{dim}(V_{g,k})$ when $\rho(g,1,k)\in \{0,1\}$ (depending on the parity of $g$).
We suspect that $V_{g,k}$ depends on very few moduli and if $g$ is suitably large we expect that $V_{g,k}=\emptyset$.
\vskip10pt
\noindent University of Michigan, Department of Mathematics\newline
East Hall, 525 East University, Ann Arbor, MI 48109-1109 \newline
e-mail: {\tt gfarkas@umich.edu }

\end{document}